\documentclass{amsart}

\usepackage{graphicx, subfigure}
\usepackage{subfloat}
\newtheorem{theorem}{Theorem}[section]

\theoremstyle{definition}

\newtheorem{example}[theorem]{Example}

\theoremstyle{remark}

\numberwithin{equation}{section}
\usepackage{cite}

\usepackage{color}
\usepackage{latexsym}
\usepackage{amsmath}
\usepackage{amssymb}
\usepackage{graphicx}
\graphicspath{%
    {converted_graphics/}
    {C:/Users/Kurdistan/Pictures/}
    {/}
}

     \title[The Computation of Zeros of Ahlfors Map... ]{The Computation of Zeros of Ahlfors Map for Doubly Connected  Regions}

     \author{Ali W.K. Sangawi}
\address{UTM Centre for Industrial and Applied Mathematics(UTM-CIAM), Universiti Teknologi Malaysia, 81310 UTM Johor Bahru, Johor, Malaysia.}
\address{Department of Mathematics, School of Science, Faculty of Science and Science Education, University of Sulaimani, 46001 Sulaimani, Kurdistan, Iraq.}
\email{alisangawi2000@yahoo.com (Ali W. K. Sangawi)} 
\author{Kashif Nazar}
     \address{Department of Mathematical Sciences, Faculty of Science,
Universiti Teknologi Malaysia, 81310 UTM Johor Bahru, Johor, Malaysia.}
\address{Department of Mathematics, COMSATS Institute of Information Technology,
P.O.Box 54000 Defence Road Off Raiwind Road Lahore Pakistan.}
\email{knazar@ciitlahore.edu.pk (Kashif Nazar)}    
\author{Ali H.M. Murid $^*$}
\address{UTM Centre for Industrial and Applied Mathematics(UTM-CIAM), Universiti Teknologi Malaysia, 81310 UTM Johor Bahru, Johor, Malaysia.}
\address{Department of Mathematical Sciences, Faculty of Science,
Universiti Teknologi Malaysia, 81310 UTM Johor Bahru, Johor, Malaysia.}
\email{alihassan@utm.my (Ali H. M. Murid)}

\thanks{$^*$Corresponding author at:{\ }{UTM Centre for Industrial and Applied Mathematics(UTM-CIAM), Universiti Teknologi Malaysia, 81310 UTM Johor Bahru, Johor, Malaysia.}}

\subjclass{30C30; 30E25; 65E05}
\keywords{Ahlfors map; Adjoint Neumann kernel; Generalized Neumann kernel.}  

\begin{document}
     
     \begin{abstract}
   The relation between the Ahlfors map and Szeg\"o kernel $S(z, a)$ is classical. The Szeg\"o kernel is a solution of a Fredholm integral equation of the second kind with the Kerzman-Stein kernel. The exact zeros of the Ahlfors map are unknown except for the annulus region. This paper presents a numerical method for computing the zeros of the Ahlfors map of any bounded doubly connected region. The method depends on the values of $S(z(t),a)$, $S'(z(t),a)$ and $\theta'(t)$ where  $\theta(t)$ is the boundary correspondence function of Ahlfors map. A formula is derived for computing $S'(z(t),a)$. An integral equation is constructed for solving $\theta'(t)$. The numerical examples presented here prove the effectiveness of the proposed method.
     \end{abstract}
     \maketitle
\section{Introduction}

The conformal mapping from a multiply connected region onto  the unit disk is known as the  Ahlfors map. If the region is simply connected then the Ahlfors map reduces to the Riemann map. Many of the geometrical features of a Riemann mapping function are shared with Ahlfors map.
The Riemann mapping function can be regarded as a solution of the following extremal problem:
For a simply connected region $\Omega$ and canonical region $D$ in the complex plane $\mathbb{C}$ and fixed $a$ in $\Omega$, construct an extremal analytic map
\begin{equation*}\label{eq1}
F:\Omega \rightarrow D ~ {\rm with} \quad F'(a)>0.
\end{equation*}
The Riemann map is the solution of this problem. It is unique conformal, one-to-one and onto map with $F(a)=0$.

   For a  multiply connected region $\Omega$ of connectivity $n>1$, the answer to the same extremal problem above becomes the Ahlfors map. It is unique analytic map
\begin{equation*}
f : \Omega \rightarrow D
\end{equation*}
that is onto, $f'(a)>0$ and $f(a)=0$.
However it has $2n-2$ branch points in the interior and is no longer one-to-one there. In fact it maps $\Omega$ onto $D$ in an \textit{n}-to-one fashion, and maps each boundary curve one-to-one onto the unit circle (see ~\cite{kra} and ~\cite{teg}). Therefore the Ahlfors map can be regarded as the Riemann mapping function in the multiply connected region.

 Conformal mapping of multiply connected regions can be computed efficiently using  the integral equation method. The integral equation method has been used by many authors to compute the one-to-one conformal mapping from multiply connected regions onto some standard canonical regions ~\cite{ker78,ker86,lee,nas-fun,nas-siam,nas-jmaa,nas-jmaa2,nas-siam2,odo,sangawi12,Sangawi11,sangawi2012circular,SangawiIJSER,SangawiAMC,arif}.

Some integral equations for Ahlfors map have been given in ~\cite{bel86,mur,Nas13,teg98,teg}. In ~\cite{ker78}, Kerzman and Stein have  derived a uniquely solvable boundary integral equation  for computing the Szeg\"o kernel of a bounded region and this method has been generalized in ~\cite{bel86} to compute Ahlfors map of bounded multiply connected regions without relying on the zeros of Ahlfors map.
In ~\cite{mur,Nas13} the  integral equations for Ahlfors map of doubly connected regions requires knowledge of zeros of Ahlfors map, which are unknown in general.

In this paper, we extend the approach of Sangawi ~\cite{sangawi12,SangawiIJSER,SangawiAMC,san-st} to construct an integral equation  for solving $\theta'(t)$ where  $\theta(t)$ is the boundary correspondence function of Ahlfors map of multiply connected region onto a unit disk. 

The plan of this paper is as follows: Section~\ref{Aux} presents some auxiliary materials. We shall derive in Section~\ref{IETht}, a boundary integral equation satisfied by $\theta'(t)$, where $\theta(t)$ is the boundary correspondence function of  Ahlfors map of bounded multiply connected regions onto a disk. From the computed values of $\theta'(t)$, we then determine the Ahlfors map. In Section~\ref{IESp}, we present a formula for computing the derivative of Szeg\"o kernel which we used for finding the derivative of Ahlfors map giving another way of computing $\theta'(t)$ analytically. In Section~\ref{2ndZero}, we present a method for computing  the second zero of the Ahlfors map for doubly connected regions.  In Section~\ref{Num}, we present some examples to illustrate our boundary integral equation method for finding the zeros of Ahlfors map for general doubly connected region. The numerical examples are first  restricted to annulus region for which the exact Ahlfors map is known, then  verified on general doubly connected region and obtained accurate results. Finally, Section~\ref{Concls} presents a short conclusion. 
\section{Auxiliary Materials}
\label{Aux}
\begin{figure}[h!] %
\centerline{\scalebox{0.4}[0.4]{\includegraphics{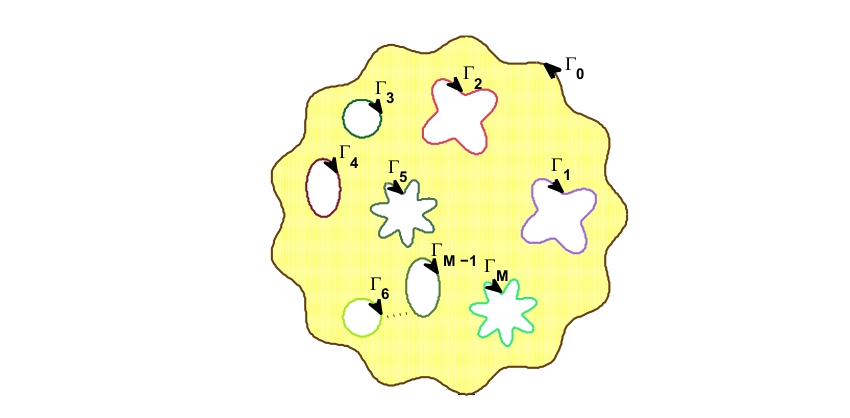}}}
\caption{A bounded multiply connected region $\Omega$ of connectivity ${M}+1$.} 
\label{fig1}
\end{figure}

Let $\Omega $ be a bounded multiply connected region of connectivity $\textit{M}+1$. The boundary $\Gamma$ consists of $M+1$ smooth Jordan curves $\Gamma_0,\Gamma_1,\dots,\Gamma_{{M}}$ such that $\Gamma_1,\dots ,\Gamma_{M}$ lie in the interior of $\Gamma_0$, where the outer curve $\Gamma_0$ has counterclockwise orientation and inner curves $\Gamma_1,\dots ,\Gamma_{\textit{M}}$  have clockwise orientation. The positive direction of the contour $\displaystyle \Gamma= \Gamma _{0}\cup \Gamma _{1}\cup \dots \cup  \Gamma _{M}$ is usually that for which $\Omega$ is on the left as one traces the boundary as shown in Figure~\ref{fig1}.

The curves $\Gamma_j$ are parameterized by $2\pi$-periodic triple continuously differentiable complex-valued functions $z_j(t)$  
with non-vanishing first derivatives
\begin{eqnarray*}
z'_j(t)=dz_j(t)/dt\not=0,\quad  t\in J_j=[0,2\pi],\quad j=0,1 ,\dots ,M. 
\end{eqnarray*}
The total parameter domain $J$ is defined as the disjoint union of $M+1$ intervals $J_0 ,\dots, J_M$. The notation
\begin{eqnarray*}\label{eeqq1}
 z(t)=z_j(t), \quad t \in J_j, \quad j=0,1, \dots ,M.
\end{eqnarray*}
must be interpreted as follows~\cite{Nas13}: For a given $\tilde t\in[0,2\pi]$, to evaluate the value of $z(t)$ at $\tilde t$, we should know in advance the interval $J_j$ to which $\tilde t$ belongs, i.e., we should know the boundary $\Gamma_j$ contains $z(\tilde t)$, then we compute $z(\tilde t)=z_j(\tilde t)$.

The generalized Neumann kernel formed with a complex continuously differentiable $2\pi$-periodic function $A(t)$  for all $t\in J$, is defined by ~\cite{WagNas,wegm}
\begin{equation*}
N(t,s)=\frac{1}{\pi}{\textrm{Im}} \left(\displaystyle \frac{A(t)}{A(s)}\frac{z'(s)}{z(s)-z(t)}\right).
\end{equation*}
The classical Neumann kernel is the generalized Neumann kernel formed with $A(t)=1$, i.e. 
\begin{eqnarray*}
N(t,s)=\frac{1}{\pi}\textrm{Im} \left(\frac{z'(s)}{z(s)-z(t)} \right).
\end{eqnarray*}
The kernel is continuous which takes on the diagonal the values 
\begin{equation*}
N(t,t)=\frac{1}{\pi} \left(\displaystyle \frac{1}{2}\textrm{Im}\frac{z''(t)}{z'(t)}-\textrm{Im}\frac{A'(t)}{A(t)}\right).
\end{equation*}

The adjoint of the function $A$ is given by 
\begin{eqnarray*}
\tilde{A}=\frac{z'(t)}{A(t)}.
\end{eqnarray*}
The generalized Neumann kernel $\tilde{N}(s,t)$ formed with $\tilde{A}$ is given by 
\begin{equation*}
\tilde{N}(s,t)=\frac{1}{\pi}\textrm{Im} \left(\displaystyle \frac{\tilde{A}(s)}{\tilde{A}(t)}\frac{z'(t)}{z(t)-z(s)}\right),
\end{equation*}
which implies
\begin{equation*}
\tilde{N}(s,t)=-N^{*}(s,t),
\end{equation*}
where $N^{*}(s,t)=N(t,s)$ is the adjoint kernel for the generalized Neumann kernel $N(s,t)$ (see~\cite{WagNas}). 

In the remaining part of this paper, we shall only be dealing with $A(t)=1$, which the classical Neumann kernel. The adjoint kernel $N^{*}(s, t)$ of the classical Neumann kernel is given by 
\begin{eqnarray*}
N^{*}(t,s)=N(s,t)=\frac{1}{\pi}\textrm{Im} \left(\frac{z'(t)}{z(t)-z(s)} \right).
\end{eqnarray*}
It is known that $\lambda =1$ is an eigenvalue of the kernel $N$ with multiplicity $1$ and $\lambda =-1$ is an eigenvalue of the kernel $N$ with multiplicity $M$ \cite{wegm}. The eigenfunctions of $N$ corresponding to the eigenvalue $\lambda = -1$ are $\left\{ \chi^{[1]}, \chi^{[2]}, \dots, \chi^{[M]} \right\}$, where
\begin{eqnarray*}
\chi^{[j]}(\xi)=\left\{ \begin{array}{cl}
          \displaystyle 1, & \displaystyle \hspace{.2cm} \xi \in \Gamma_{j},\\
          \displaystyle  0, & \displaystyle \hspace{.2cm} \rm{otherwise},
            \end{array}
       \right.
\hspace{.5cm}j=1, 2, \dots, M.
\end{eqnarray*}
Let $H^{*}$ be the space of all real H\"older continuous $2\pi$- periodic functions $\omega(t)$ of the parameter $t$ on $J_j$ for $j=0,1,\dots,M$, $\textrm{i.e}.$
\begin{equation*}
\omega(t)=\omega_k(t),~ t\in J_j,~k=0,1,\dots,M, j=0,1,\dots,M.
\end{equation*}
Define also the integral operator $\textbf{J}$ by \cite{Sangawi11,SangawiIJSER} 
\begin{equation*}
\textbf{J} {\nu}= \int_J \frac{1}{2\pi}\sum_{j=1}^M \chi^{[j]}(s)\chi^{[j]}(t)\nu(s)ds.
\end{equation*}
We also define the Fredholm integral operator $N^*$ by
\begin{equation*}
\textbf{N}^*\psi(t)=\int_J N^*(t,s)\psi(s)ds,~ t\in J.
\end{equation*}
A complex-valued function $P(z)$ is said to satisfy the interior relationship if $P(z)$  is analytic in $\Omega$ and satisfies the non-homogeneous boundary relationship
\begin{equation}\label{eq6a}
P(z)=\frac{b(z) \overline{T(z)}}{\overline{G(z)}}\overline{P(z)}+\overline{H(z)}, \quad z\in \Gamma, 
\end{equation}
where $G(z)$ is analytic in $\Omega$ and H\"older continuous on $\Gamma$, $G(z)\not=0$ on $\Gamma$, and $T(z) = z'(t)/\left|z'(t)\right|$ denotes the unit tangent in the direction of increasing parameters at the point $z=z(t)\in\Gamma$. The boundary relationship (\ref{eq6a}) also has the following equivalent form:
\begin{equation}\label{eq7}
G(z)=\overline{b(z)}T(z)\frac{P(z)^2}{|P(z)|^2}+\frac{G(z)H(z)}{\overline{P(z)}}, \quad z\in \Gamma.
\end{equation}

The following theorem gives an integral equation for an analytic function satisfying the non-homogeneous boundary relationship (\ref{eq6a}) or (\ref{eq7}).\\
\begin{theorem}\label{theo2}
\cite{sangawi12} If the function $P(z)$ satisfies the interior relationship with (\ref{eq6a}) or (\ref{eq7}), then
\begin{eqnarray*}\label{eq8}
&&{}T(z)P(z)+\int_{\Gamma} K(z,w)T(w)P(w)|dw| \nonumber\\
&&{}\qquad\qquad+ b(z)\left[\displaystyle\sum_{a_j {\textrm{inside}} \Gamma} \underset{w=a_j}{\textrm{Res}}  \frac{P(w)}{(w-z)G(w)}\right]^- =-\overline{L^{-}_{R}(z)},\nonumber\\
\end{eqnarray*}
where
\begin{equation*}\label{eq9}
K(z,w)=\frac{1}{2\pi {\rm{i}}}\left[\displaystyle\frac{T(z)}{z-w}- \frac{b(z)\overline{T(w)}}{b(w)(\overline{z}-\overline{w})}\right],
\end{equation*}
and 
\begin{equation*}\label{eq10}
L^{-}_{R}(z)=\frac{-1}{2}\frac{H(z)}{T(z)}+{\rm PV} \frac{1}{2\pi{\rm{i}}}\int_{\Gamma}\frac{\overline{b(z)}H(w)}{\overline{b(w)}(w-z)T(w)}dw.
\end{equation*}
The symbol $``-"$ in the superscript denotes the complex conjugate and the sum is over all those zeros  that lie inside $\Omega$. If $G$ has no zeros in $\Omega$, then the term containing the residue will not appear.
\end{theorem}
\section{Integral equation method for computing $\theta'$}\label{IETht}
Let $f(z)$ be the Ahlfors 
function which maps $\Omega$ conformally  onto a unit disc.
The mapping function $f$ is determined up to a factor of modulus 1. 
The function $f$ could be made unique by prescribing that
\begin{eqnarray*}
f(a_j) = 0,\quad f'(a_0) > 0, \quad  j=0,1,2,\dots ,M. 
\end{eqnarray*}
where $a_j \in\Omega ,j=0,1,\dots, M$ are the zeros of the Ahlfors map. The boundary values of $f$ can be represented as
\begin{eqnarray} \label{eq1}
f(z_j(t))=e^{{\rm{i}}{\rm \theta} _{j}(t)},\quad \Gamma _j:z=z_j(t),\quad 0 \le t \le \beta _{j},
\end{eqnarray}
where $\theta _{j}(t), j=0,1 ,\dots ,M $ are the boundary correspondence functions of
$\Gamma _{j}$.
Also we have  
\begin{eqnarray}\label{eq27a}
\frac{f'(z(t))z'(t)}{f(z(t))}={\rm i}\theta'(t),
\end{eqnarray}

The unit tangent to $\Gamma$ at $z(t)$ is denoted by $T(z(t)) = z'(t)/|z'(t)|$. Thus it can be shown that
\begin{equation}\label{eq3}
f(z_{j}(t))=\frac{1}{{{\rm i}}} T(z_{j}(t))\frac{\theta _{j}'(t)}{|\theta _{j}'(t)|}\frac{f'(z_{j}(t))}{|f'(z_{j}(t))|},\qquad z_j\in\Gamma_j.
\end{equation}
By the angle preserving property of conformal  map, the image of $\Gamma_0$ remains in counter-clockwise orientation so $\theta _{0}'(t) > 0$ \
while the images of inner boundaries $\Gamma_j$  in clockwise orientation so $\theta' _{j}(t)< 0$, for $j= 1,\dots ,M$.
Thus\\
\begin{equation*}\label{e:piece}
\textrm{sign}(\theta'_j(t)) = \left\{
\begin{array}{l@{\hspace{0.5cm}}l}
+1,     &j=0, \\
-1,     &j=1,\dots ,M. \\
\end{array}%
\right.
\end{equation*}
The boundary relationship (\ref{eq3}) can be written briefly as
\begin{equation} \label{eq5}
f(z)=\textrm{sign}(\theta '(t))\frac{1}{{{\rm i}}} T(z)\frac{f'(z)}{|f'(z)|} , \quad z\in\Gamma .
\end{equation}  
Since the Ahlfors map can be written as $\displaystyle f(z)= \prod_{j=0}^{M}(z-a_j) \hat g(z)$,  where $\hat g(z)$ is analytic in $\Omega$ and $\hat g(z)\not=0$ in $\Omega$, 
 we have 
\begin{eqnarray} \label{eq11}
\frac{f'(z)}{f(z)}=\frac{\hat g'(z)}{\hat g(z)}+\sum_{j=0}^{M}\frac{1}{z-a_j}, \quad z\in\Gamma ,  
\end{eqnarray}
so that
\begin{eqnarray}\label{eq12}
 D(z)=\frac{f'(z)}{f(z)}-\sum_{j=0}^{M}\frac{1}{z-a_j}=\frac{\hat g'(z)}{\hat g(z)}.
\end{eqnarray}
is analytic in $\Omega$. Squaring both sides of (\ref{eq5}), gives

\begin{eqnarray} \label{eq13}
f(z)^2= -T(z)^2 \frac{f'(z)^2}{|f'(z)|^2},   
\end{eqnarray}
which implies
\begin{equation} \label{eq14}
\frac{f'(z)}{f(z)}= -\overline{T(z)}^2 \overline{ \left(\displaystyle \frac{f'(z)}{f(z)}\right)}. 
\end{equation}
From (\ref{eq12}), we have
\begin{equation} \label{eq15}
\frac{f'(z)}{f(z)}= D(z)+\sum_{j=0}^{M}\frac{1}{z-a_j}.
\end{equation}
Substituting (\ref{eq15}) into (\ref{eq14})
\begin{eqnarray}\label{eq16}
 D(z)+\sum_{j=0}^{M}\frac{1}{z-a_j}=-\overline{T(z)}^2 \overline{D(z)}-\sum_{j=0}^{M}\frac{\overline{T(z)}^2}{\overline{z}-\overline{a_j}} 
\end{eqnarray}
or
\begin{equation} \label{eq17}
 D(z)=-\overline{T(z(t))}^2 \overline{D(z)}-\sum_{j=0}^{M}\frac{\overline{T(z)^2}}{\overline{z}-\overline{a_j}}-\sum_{j=0}^{M}\frac{1}{z-a_j}.
\end{equation}
Comparing (\ref{eq17}) with (\ref{eq6a}), we get

\begin{equation}\label{eq18}
P(z)=D(z), b(z)=-\overline{T(z)}, G(z)=1,H(z)=-\sum_{j=0}^{M}\frac{T(z)^2}{z-a_j}-\sum_{j=0}^{M}\frac{1}{\overline{z}-\overline{a_j}}.
\end{equation}
Substituting these assignments into  Theorem~\ref{theo2}, we obtain
\begin{eqnarray}\label{eq19}
T(z)D(z)+\int_{\Gamma} K(z,w)T(w)D(w)|dw|
+0 =-\overline{L^{-}_{R}(z)},
\end{eqnarray}
where
\begin{equation*}\label{eq20}
\begin{array}{lll}
K(z,w)&\displaystyle=\frac{1}{2\pi {\rm i}}\left[\displaystyle\frac{T(z)}{z-w}- \frac{\overline{T(z)}}{\overline{z}-\overline{w}}\right]\vspace{6pt}
&\displaystyle=\frac{1}{\pi}{\rm{Im}} \left(\displaystyle \frac{\textit{T(z)}}{\textit{z}-\textit{w}}\right)=N(z,w)
\end{array}
\end{equation*}
and
\begin{eqnarray*}\label{eq21}
&&{}L^{-}_{R}(z)=\frac{-1}{2T(z)} \left[\displaystyle \sum_{j=0}^M\frac{T(z)^2}{z-a_j}-\sum_{j=0}^M \frac{1}{\overline{z}-\overline{a_j}}\right] \nonumber\\
&&{}\qquad\qquad+{\rm PV}\frac{1}{2\pi {\rm i}}\int_{\Gamma}\frac{T(z)}{(w-z)T(w)^2}\left[\displaystyle \sum_{j=0}^M\frac{T(w)^2}{w-a_j}-\sum_{j=0}^M \frac{1}{\overline{w}-\overline{a_j}}\right].\nonumber\\
\end{eqnarray*}
 Using the fact that \cite{Sangawi11}
\begin{eqnarray*}\label{eq22}
{\rm PV}\frac{1}{2\pi {\rm i}} \int_{\Gamma}\frac{1}{(w-z)(w-a_j)}dw=-\frac{1}{2(z-a_j)},
\end{eqnarray*}
 $dw=T(w)|dw|$, and after some simplifications, $(\ref{eq19})$ becomes
\begin{eqnarray}\label{eq23}
T(z) \frac{f'(z)}{f(z)}+\int_{\Gamma} N(z,w)\frac{f'(w)}{f(w)}T(w)|dw|=2 {\rm i} {\rm{Im}}\left[ \sum_{\textit{j}=0}^{\textit{M}} \frac{\textit{T(z)}}{\textit{z}(t)-\textit{a}_j}\right],\nonumber\\
\end{eqnarray}
where
\begin{equation}\label{eq24}
N(z,w) = \left\{
\begin{array}{l@{\hspace{0.5cm}}l}
\displaystyle \frac{1}{2\pi {\rm i}}\left[\displaystyle \frac{T(z)}{z-w}-\frac{\overline{T(z)}}{\overline{z}-\overline{w}}\right],     &\displaystyle z\not= w\in \Gamma, \\
\displaystyle \frac{1}{2\pi |z'(t)|} {\rm{Im}}\left(\displaystyle \frac{\textit{z}''(t)}{\textit{z}'(t)}\right),     &\displaystyle z=w\in \Gamma ,\\
\end{array}%
\right.
\end{equation}
with
\begin{eqnarray}\label{eq25}
\frac{1}{2\pi {\rm i}}\int_{-\Gamma_{j}}\frac{f'(w)T(w)}{f(w)}|dw|=1,\quad j=1,2,\dots ,M.
\end{eqnarray}
In  integral equation $(\ref{eq23})$, letting $\displaystyle z=z(t)$, $w=z(s)$ and multiplying both sides by $|z'(t)|$, gives

\begin{equation}\label{eq26}
\frac{f'(z(t))}{f(z(t))} z'(t) +\int_{\Gamma} N(z(t),z(s))\frac{f'(z(s))}{f(z(s))}z'(s)ds=2 {\rm i} {\rm{Im}}\left[ \sum_{j=0}^{M} \frac{\textit{z}'(t)}{\textit{z}(t)-\textit{a}_j}\right].
\end{equation}
Using (\ref{eq27a}), the above equation becomes 
\begin{eqnarray}\label{eq28}
\theta'(t)+\int_{J}N(t,s)\theta'(s)ds=2{\rm Im}\left[\sum_{j=0}^{M}\frac{\textit{z}'(t)}{\textit{z}(t)-\textit{a}_j}\right].
\end{eqnarray}
Since $\displaystyle N^{*}(s,t)=N(t,s)$, the integral equation $(\ref{eq28})$ in the operator form is 
\begin{eqnarray}\label{eq29}
(\textbf{I}+\textbf{N}^{*})\theta'_j = \phi , \qquad j=0,1,\dots ,M,
\end{eqnarray}
where
\begin{eqnarray}\label{eq30}
\phi= 2{\rm Im} \left[\displaystyle \sum_{j=0}^{M}\frac{\textit{z}'(t)}{\textit{z}(t)-a_j}\right].
\end{eqnarray}
By Theorem 12 in ~\cite{wegm},  $\lambda=-1$ is an eigenvalue of $N^*$ with multiplicity $M$, therefore the integral equation (\ref{eq29}) is not solvable. To overcome this problem, we note from (\ref{eq25}) and (\ref{eq27a}) that

\begin{equation*}\label{eq31}
\textbf{J} \theta'=\left(\displaystyle 0, \frac{1}{2\pi}\int_{J_{1}}\theta'(s)ds, \frac{1}{2\pi}\int_{J_{2}}\theta'(s)ds,... , \frac{1}{2\pi}\int_{J_{M}}\theta'(s)ds\right)
\end{equation*}
which implies  that 
\begin{eqnarray}\label{eq32}
\textbf{J} \theta'=\psi,
\end{eqnarray}
where
\begin{eqnarray}\label{eq33}
\psi = (0,1,1,\dots,1).
\end{eqnarray}
By adding (\ref{eq29}) and (\ref{eq32}), we get
\begin{eqnarray}\label{eq34}
(\textbf{I}+\textbf{N}^*+\textbf{J})\theta'=\phi+\psi 
\end{eqnarray}
which implies
\begin{equation}\label{phio}
\phi=(\textbf{I}+\textbf{N}^*+\textbf{J})\theta'-\psi.
\end{equation}

In the next section, we derive another formula for computing $\theta'$ from which the second zero of the Ahlfors map for doubly connected region can be calculated using (\ref{phio}).
\section{An Alternative Formula for Computing $\theta'$} \label{IESp}
  The Ahlfors map is related to the Szeg\"o kernel $S(z,a_0)$ and the Garabedian kernel $L(z,a_0)$ by \cite{bel86}
\begin{equation}\label{belAhlfora}
f(z)=\frac{S(z,a_0)}{L(z,a_0)}, \quad z\in \Omega \cup \Gamma.
\end{equation}
The Szeg\"o kernel $S(z,a_0)$ and Garabedian kernel $L(z,a_0)$ are related on $\Gamma$ as
\begin{equation*}
L(z,a_0)=-{\rm i} \overline{T(z)}\overline{S(z,a_0)},\quad z\in \Gamma,
\end{equation*}
so that (\ref{belAhlfora}) becomes
\begin{equation}\label{Ahlsimpl}
f(z)=\frac{1}{{\rm i}}\frac{S(z,a_0)T(z)}{\overline{S(z,a_0)}},\quad z\in \Gamma,
\end{equation}
or
\begin{equation}\label{Ahlforsfoft}
f(z(t))=\frac{1}{{\rm i}}\frac{S(z(t),a_0)T(z(t))}{\overline{S(z(t),a_0)}},\quad z(t)\in\Gamma.
\end{equation}
From Bell \cite{bel86}, we have the integral equation for the Szeg\"o kernel
\begin{equation}\label{eq11a}
S(z,a_0)+\int_{\Gamma} A(z,w)S(w,a_0)|dw|=g(z),\quad z\in \Gamma,
\end{equation}
where
\begin{eqnarray*}
A(z,w)=\left\{\begin{array}{ll}
\displaystyle\frac{1}{2\pi{\rm i}}\displaystyle\left[\frac{T(w)}{z-w}-\frac{\overline{T(z)}}{\overline{z-w}}\right],&\displaystyle \quad z\not=w \in\Gamma\\
\displaystyle0,&\displaystyle \quad z=w, 
\end{array}\right.               
\end{eqnarray*}
and 
\begin{equation}\label{eq34a}
g(z(t))=-\frac{1}{2\pi {\rm i}}\frac{\overline{T(z(t))}}{\overline{z(t)-a_0}}.
\end{equation}
With $z=z(t)$ and $w=z(s)$, (\ref{eq11a}) becomes
\begin{equation*}
S(z(t),a_0)+\int_{J_j} A(z(t),z(s))S(z(s),a_0)|z'(s)|ds=g(z(t)).
\end{equation*}
Differentiate both sides with respect to $t$, we get
\begin{equation*}
\frac{d}{dt}S(z(t),a_0)+\frac{d}{d t}\int_{J_j}A(z(t),z(s)) S(z(s),a_0)|z'(s)|ds=\frac{d}{dt}g(z(t)),
\end{equation*}
which is equivalent to
\begin{equation}\label{sp}
S'(z(t),a_0)z'(t)=g'(z(t))z'(t)-\int_{J_j}\displaystyle\left[\frac{d}{d t}A(z(t),z(s))z'(t)\right] S(z(s),a_0)|z'(s)|ds.
\end{equation}
Differentiate (\ref{eq34a}) with respect to $t$, gives
\begin{equation}\label{gp}
g'(z(t))z'(t)=-\frac{1}{2\pi {\rm i}}\displaystyle\left[\frac{\overline{T'(z(t))z'(t)}}{\overline{z(t)-a_0}}-\frac{\overline{z'(t)}~\overline{T(z(t))}}{\overline{(z(t)-a_0)^2}}\right].
\end{equation}
Since
\begin{equation*}
T(z(t))=\frac{z'(t)}{|z'(t)|},
\end{equation*}
we get
\begin{equation}
T'(z(t))z'(t)=\frac{z''(t)}{2|z'(t)|}-\frac{z'^2(t)\overline{z''(t)}}{2|z'(t)|^3}.
\end{equation}
Since 
\begin{equation*}
A(z(t),z(s))=\frac{1}{2\pi {\rm i}}\displaystyle\left[\frac{T(z(s))}{z(t)-z(s)}-\frac{\overline{T(z(t))}}{\overline{z(t)-z(s)}}\right],\quad {\rm for}~ t\neq s,
\end{equation*}
differentiating with respect to $t$, we get
\begin{eqnarray*}
\frac{d}{dt}A(z(t),z(s))
=\frac{1}{2\pi {\rm {\rm i}}}\displaystyle \left[\frac{-z'(t)~T(z(s))}{(z(t)-z(s))^2}-\frac{\overline{T'(z(t))z'(t)}}{(\overline{z(t)}-\overline{z(s)})}+\frac{\overline{T(z(t))}~\overline{z'(t)}}{(\overline{z(t)}-\overline{z(s)})^2}\right],
\end{eqnarray*}
which implies
\begin{equation}\label{eqqker}
\begin{array}{llll}
\displaystyle 2\pi {\rm i}|z'(s)|\frac{d}{dt}A(z(t),z(s))=&-&\displaystyle\left(\frac{z'(t)}{z(t)-z(s)}\right)^2 \displaystyle\left(\frac{z'(s)}{z'(t)}\right)\\
&-&\displaystyle\left(\frac{\overline{z''(t)}z'(t)-\overline{z'(t)}z''(t)}{2|z'(t)|^2 }\right) \displaystyle\left(\overline{\frac{z'(t)}{z(t)-z(s)}}\right)\displaystyle\left|\frac{z'(s)}{z'(t)}\right|\\
&+&\displaystyle\left(\overline{\frac{z'(t)}{z(t)-z(s)}}\right)^2\displaystyle\left|\frac{z'(s)}{z'(t)}\right|.
\end{array}
\end{equation}
For  $s$ sufficiently close to $t$, we have
\begin{equation*}
z(s)=z(t)+z'(t)(s-t)+\frac{1}{2}z''(t)(s-t)^2+\frac{1}{6}z'''(t)(s-t)^3+O((s-t)^4)
\end{equation*}
which implies
\begin{equation}\label{eqq1ker}
\frac{z'(t)}{z(s)-z(t)}=\frac{1}{s-t}\displaystyle\left[1+\frac{1}{2}\frac{z''(t)}{z'(t)}(s-t)+\frac{1}{6}\frac{z'''(t)}{z'(t)}(s-t)^2+O((s-t)^4)\right]^{-1}
\end{equation}
Using the fact that 
\begin{eqnarray}\label{eqq2ker}
\frac{1}{1+\varrho}=1-\varrho+\varrho^2+O(\varrho^3),
\end{eqnarray}
for $\varrho$ close to zero, (\ref{eqq1ker}) becomes 
\begin{equation}\label{eqq3ker}
\frac{z'(t)}{z(s)-z(t)}=\frac{1}{s-t}-\frac{1}{2}\frac{z''(t)}{z'(t)}-\frac{1}{6}\frac{z'''(t)}{z'(t)}(s-t)+\frac{1}{4}\frac{z''^2(t)}{z'^2(t)}(s-t)+O((s-t)^2).
\end{equation}
We observe that 
\begin{equation}\label{eqq4ker}
\frac{z'(s)}{z'(t)}=1+\frac{z''(t)}{z'(t)}(s-t)+\frac{1}{2}\frac{z'''(t)}{z'(t)}(s-t)^2+O((s-t)^3)
\end{equation}
Next observe that
\begin{equation}\label{eqq5ker}
\begin{array}{lll}
\displaystyle\left|\frac{z'(s)}{z'(t)}\right|^2=\displaystyle 1&\displaystyle +&\displaystyle 2{\rm Re}\left(\frac{z''(t)}{z'(t)}\right)(s-t)\\
&\displaystyle +&\displaystyle \left({\rm Re}\left(\frac{z'''(t)}{z'(t)}\right)+\left|\frac{z''(t)}{z'(t)}\right|^2\right)(s-t)^2+O((s-t)^3)
\end{array}
\end{equation}
Using the fact that 
\begin{equation}\label{eqq6ker}
\sqrt{1+\hat x}=1+\frac{1}{2}\hat x-\frac{1}{8}\hat x^2+\cdots
\end{equation}
for small $\hat x$, (\ref{eqq5ker}) becomes
\begin{equation}\label{eqq7ker}
\displaystyle\left|\frac{z'(s)}{z'(t)}\right|=1+{\rm Re}\displaystyle\left(\frac{z''(t)}{z'(t)}\right)(s-t)+\frac{1}{2}\left({\rm Re}\displaystyle\left(\frac{z'''(t)}{z'(t)}\right)+\displaystyle\left|\frac{z''(t)}{z'(t)}\right|^2\right)(s-t)^2+O((s-t)^3).
\end{equation}
Substituting (\ref{eqq7ker}), (\ref{eqq4ker}), (\ref{eqq3ker}) into (\ref{eqqker}) and then taking the limit $s \to t$ to both sides, we obtain 
\begin{equation*}\label{kerA1}
\frac{d}{dt}A(z(t),z(s)) = \left\{
\begin{array}{l@{\hspace{0.5cm}}l}
\displaystyle \frac{1}{2\pi {\rm i}}\displaystyle \left[\frac{-z'(t)~T(z(s))}{(z(t)-z(s))^2}-\frac{\overline{T'(z(t))}~\overline{z'(t)}}{(\overline{z(t)}-\overline{z(s)})}+\frac{\overline{T(z(t))}~\overline{z'(t)}}{(\overline{z(t)}-\overline{z(s)})^2}\right],\vspace{4pt}\\
\hfill t\not= s\in \Gamma, \\
 \displaystyle \frac{{1}}{12 \pi \left|z'(t)\right|}{\rm Im}\displaystyle \left(\frac{z'''(t)}{z'(t)}\right),\hfill t=s\in \Gamma.\\
\end{array}%
\right.
\end{equation*}
Using the values from (\ref{eq11a})  and  (\ref{gp}), we can find the derivative of Szeg\"o kernel from (\ref{sp}).
By defining the following terms, 
\begin{eqnarray}\label{pr}
&&{} \displaystyle f_p=f'(z(t))z'(t),\quad S_p=S'(z(t),a_0)z'(t), \nonumber\\
&&{}\displaystyle  S_z=S(z(t),a_0),\quad T_p=T'(z(t))z'(t),\quad T_z=T(z(t)).
\end{eqnarray}
and differentiating both sides of (\ref{Ahlforsfoft}) with respect to $t$, the result is 
\begin {equation}\label{ahlforsp}
f_p=\frac{1}{{\rm {\rm i}}}\displaystyle\left[\frac{S_p~T_z+S_z~T_p}{\overline{S_z}}-\frac{\overline{S_p}~(S_z T_z)}{\overline{(S_z)^2}}\right].
\end{equation}
By using the values from (\ref{Ahlsimpl}) and (\ref{ahlforsp}) in (\ref{eq27a}), we get 
\begin{equation}\label{thetap}
\theta'(t)=2{\rm Im}\displaystyle\left(\frac{S_p}{S_z}\right)+{\rm Im}\displaystyle\left(\frac{z''(t)}{z(t)}\right)
\end{equation}
which is the second formula for computing $\theta'$.
\section{Computing the Second Zero of the Ahlfors Map for Doubly Connected Region}\label{2ndZero}

In particular, if $\Omega$ is a doubly connected region, i.e. $M=1$, then (\ref{eq30}) becomes

\begin{equation}\label{phi2b}
\phi=2 {\rm Im} \displaystyle \left[\frac{z'(t)}{z(t)-a_0}+\frac{z'(t)}{z(t)-a_1}\right].
\end{equation}
As $\phi$ is known from (\ref{phio}), and the zero $a_0\in \Omega$ can be freely prescribed, the only unknown  in (\ref{phi2b}) is the second zero $a_1$ of Ahlfors map. The above equation can be written as

\begin{equation}\label{phi2b2}
{\rm Im} \displaystyle \left[\frac{z'(t)}{z(t)-a_1}\right]=k_1(t).
\end{equation}
where
\begin{equation}\label{phi2bb}
k_1(t)=\frac{1}{2}\displaystyle\left[\phi-2 {\rm Im} \displaystyle \left(\frac{z'(t)}{z(t)-a_0}\right)\right].
\end{equation}
Now we suppose that
\begin{eqnarray*}
&&{}z(t)=x(t)+{\rm i}y(t),\nonumber\\
&&{}a_1= \alpha+{\rm i}\beta.
\end{eqnarray*}
Then (\ref{phi2b2}) becomes
\begin{equation}
\frac{y'(t)(x(t)-\alpha)+x'(t)(\beta-y(t))}{(x(t)-\alpha)^2+(y(t)-\beta)^2}=k_1(t)
\end{equation}
After some algebraic manipulations, we obtain
\begin{equation}\label{kkk}
k_2(t)\alpha+k_3(t)\beta+k_1(t)(\alpha^2 +\beta^2)=-k_4(t),
\end{equation}
where
\begin{eqnarray}\label{k2}
k_2(t)&=&y'(t)-2k_1(t)x(t),
\end{eqnarray}
\begin{eqnarray}\label{k3}
k_3(t)&=&-2k_1(t)y(t)-x'(t),
\end{eqnarray}
\begin{eqnarray}\label{k4}
k_4(t)&=&k_1(t)x^2+k_1(t)y^2+x'(t)y(t)-y'(t)x(t).
\end{eqnarray}

Equation (\ref{kkk}) can be viewed as containing three unknowns namely $\alpha$, $\beta$ and $\alpha^2+\beta^2$. To determine them, we need three equations. For this we can choose any three values of parameter $t$ in the given interval. At these three random values of ${\it t}$, we can find the values of the coefficients $k_1(t)$, $k_2(t)$, $k_3(t)$ and $k_4(t)$. By solving the obtained system of three equations with three unknowns, we can find the  value of $\alpha$ and $\beta$, which are the real and imaginary parts of the second zero of Ahlfors map for doubly connected region.
By solving the obtained system of three equations with three unknowns, we can find the  value of $a_1=\alpha+{\rm i}\beta$, which is the second zero of Ahlfors map for doubly connected region.
\section{Numerical Examples}\label{Num}

For solving the integral equation (\ref{eq11a}) numerically, the reliable procedure is by using the Nystr\"om method with the trapezoidal rule with \textit{n} equidistant nodes in each interval $J_j$, $j=1,\dots ,M$  ~\cite{nas-fun,nas-siam,nas-jmaa,nas-jmaa2,nas-siam2}. The trapezoidal rule is the most accurate method for integrating periodic functions numerically \cite[pp.134-142]{Davis}. By solving the integral equation (\ref{eq11a}) for $S(z(t),a_0)$ gives $S_p=S'(z(t),a_0)z'(t)$ from (\ref{sp}) and $\theta'(t)$ from (\ref{thetap}). By (\ref{phio}) we get the value of $\phi$ and then from (\ref{phi2b}) we get the  value of $a_1=\alpha+{\rm i}\beta$.
We then apply (\ref{Ahlforsfoft}) and the Cauchy integral formula to compute $f(\alpha+{\rm i}\beta)$.
 For evaluating the Cauchy integral formula $f(z)=(1/(2\pi {\rm i}))\int_{\Gamma} {f(w)}/{(w-z)}dw$ numerically, we use the equivalent form
\begin{eqnarray}\label{eq39}
f(z)=\frac{\int_\Gamma \frac{f(w)}{w-z}dw}{\int_\Gamma \frac{1}{w-z}dw} ,\quad z\in \Omega,
\end{eqnarray}
which also works very well for $z\in \Omega$ near the boundary $\Gamma$. When the trapezoidal rule is applied to the integrals in (\ref{eq39}), the term in the denominator compensates the error in the numerator (see ~\cite{SangawiIJSER}). All the computations were done using MATLAB 7.12.0.635(R2011a). 
\begin{example}\label{Ex1} 
 Consider an annulus region bounded by the two circles 
\begin{eqnarray*}
\Gamma _{0}: {\{}z(t) &=& e^{{\rm i}t}{\}},\\
\Gamma _{1}: {\{}z(t) &=&r e^{-{\rm i}t}{\}},\qquad t: 0 \le t \le 2\pi, 0<r<1. 
\end{eqnarray*}
In~\cite{teg}, Tegtmeyer and Thomas computed the Ahlfors map using Szeg\"o and Garabedian kernels for the annulus region, where the authors have used series representations of both Szeg\"o kernel and Garabedian kernel. With these representations, they found the two zeros $a_0$ and $a_1=\displaystyle \frac{-r}{\overline{a_0}}$ for Ahlfors map, where $r$ is the radius of the inner circle. They have also considered the symmetry case where the zeros are $\displaystyle a_0=\sqrt{r}$ and $\displaystyle a_1=-\sqrt{r}$. This example has also been considered in ~\cite{Nas13} where Ahlfors map was computed using a boundary integral equation related to a Riemann-Hilbert problem. Here we shall use these values of zeros of Ahlfors map for comparison with our numerical zeros with proposed method  
in the annulus $\displaystyle r < |z| < 1$.
In this example we have chosen $r=0.1$ and the first zero $a_0=0.6$. 
See Table~1 for numerical comparison of the computed second zeros $a_{1n}$ from (\ref{kkk}) and the exact second zeros of the Ahlfors map, $a_1=-r/\overline{a_0}$. 
\end{example}

\begin{figure}[ht!]
     \begin{center}
            \includegraphics[width=0.45\textwidth]{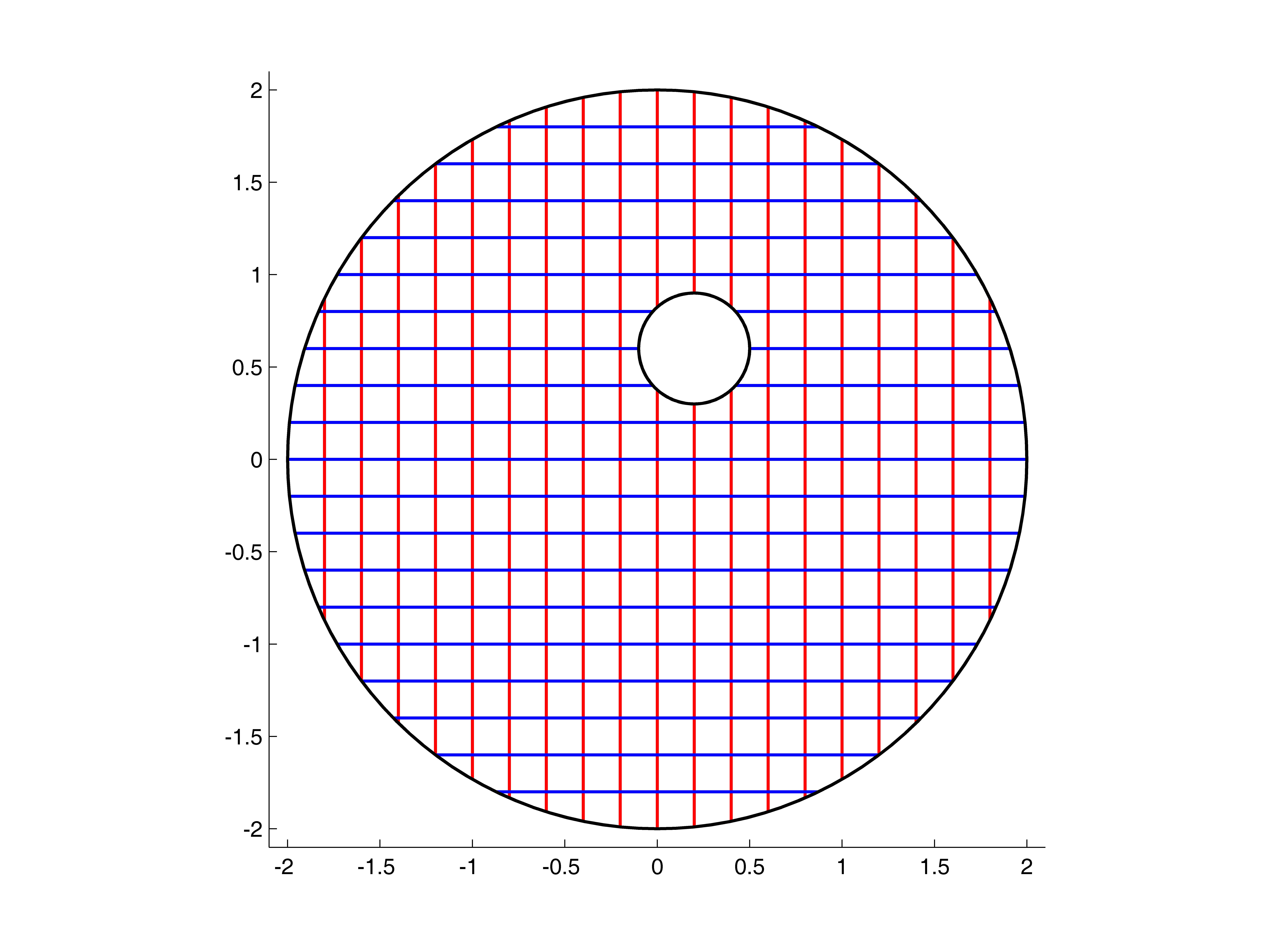}
    \end{center}
\vspace{-0.3cm}
    \caption{ The region $\Omega$ for Example~\ref{Ex2}.}%
   \label{Ex2Ahl}
\end{figure}

\begin{table}[h!]
\label{tabEx1}
\caption{Maximum error norm $\vert\vert a_{1n}-a_1\vert\vert$ with $r=0.1$ for Example~\ref{Ex1}.}
\begin{center}
\begin{tabular}{ c c c }
\hline
\multicolumn{2}{c}{\rule{0ex}{3ex}  Non-Symmetric case $(a_0=0.5, a_1=-r/\overline{a_0})$  }\\
\hline
 \textit{n}         & $\vert\vert a_{1n}-a_1\vert\vert$\\
\hline
64           &2.66(-15)\\
\hline
\hline
\multicolumn{2}{c}{\rule{0ex}{3ex}  Symmetric case   $(a_0=\sqrt{r}, a_1=-\sqrt{r})$}\\
\hline
 \textit{n}         & $\vert\vert a_{1n}-a_1\vert\vert$\\
\hline
64           &5.56(-14)\\
\hline
\end{tabular}
\end{center}
\end{table}
\begin{example}
\label{Ex2} 
\rm Consider a doubly connected region $\Omega$ bounded by the two non-concentric circles 
\begin{eqnarray*}
\Gamma _{0}: {\{}z(t) &=& 2e^{{\rm i}t}{\}},\\
\Gamma _{1}: {\{}z(t) &=& c + r e^{-{\rm i}t}{\}},\qquad t: 0 \le t \le 2\pi,
\end{eqnarray*}
with $c= 0.2+0.6{\rm i}$ and radius $r=0.3$. The test region is shown in Figure~ \ref{Ex2Ahl}.  Given a first zero $a_0$ of the Ahlfors map, the exact second zero $a_1$ is unknown for this region. We compute $a_{1n}$ from (\ref{kkk}) which is the approximate value of $a_1$. The acuracy is measured by computing $f(a_{1n})$ from (\ref{eq39}). The theoretical value of $f(a_1)$ is zero. See Table~\ref{tabEx2} for the numerical values of $a_{1n}$ and $f(a_{1n})$ involving two choices of $a_0$.
\end{example}
\begin{table}[h!]
\caption{Numerical values of $a_{1n}$ and $f(a_{1n})$ for Example~\ref{Ex2}.}\label{tabEx2}
\begin{center}
\begin{tabular}{ c c c}
\hline
$$  &$a_0=-0.5+0.2{\rm i}$ $$\\
\hline
 $n$           &$a_{1n}$            &$f(a_{1n})$\\
\hline
$16$       & $ 0.7135 + 1.0357{\rm i}$      &$2.2(-03)$\\
\hline
$32$      & $0.7125 + 1.0342{\rm i}$       &$1.12(-06)$\\
\hline
$64$      & $0.7125 + 1.0342{\rm i}$      &$2.9(-13)$\\
\hline
$128$     & $0.7125 + 1.0342{\rm i}$        &$4.8(-15)$\\
\hline
\hline
$$  &$a_0=-{\rm i}$ $$\\
\hline
 $n$           &$a_{1n}$            &$f(a_{1n})$\\
\hline
$16$      & $ 0.2831 + 1.0284{\rm i}$      &$1.0(-01)$\\
\hline
$32$      & $0.2763 + 1.0031{\rm i}$       &$7.3(-04)$\\
\hline
$64$      & $0.2763 + 1.0031{\rm i}$        &$2.7(-08)$\\
\hline
$128$     &$0.2763 + 1.0031{\rm i}$	    &$4.6(-15)$\\
\hline
\end{tabular}
\end{center}
\end{table}

\begin{figure}[ht!]
     \begin{center}
\hspace{-0.9cm}
        \subfigure[ $\Omega$~$(\xi=0.6452-0.8655{\rm i})$ ]{%
            \label{fig3:first}
            \includegraphics[width=0.45\textwidth]{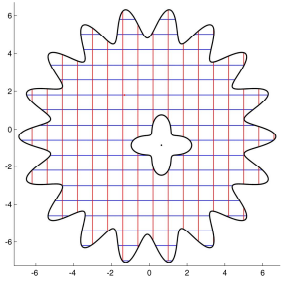}\hspace{-0.8cm}
        }%
        \subfigure[ $\Omega$~$(\xi=-2.4516+2.3626{\rm i})$ ]{%
            \label{fig3:firstt}
            \includegraphics[width=0.45\textwidth]{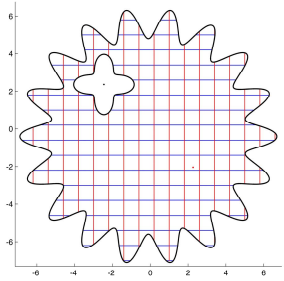}\hspace{-0.9cm}
        }%
    \end{center}
\vspace{-0.3cm}
    \caption{ The region $\Omega$ for Example~\ref{Ex3}.     }%
   \label{Ex3Ahl}
\end{figure}

\newpage
\begin{example}
\label{Ex3} 
\rm Consider a doubly connected region $\Omega$ bounded by $\Gamma_{0}$ and $\Gamma_{1}$
\begin{eqnarray*}
\Gamma _{0}: {\{}z(t) &=& -0.1-0.4i+ (6+0.8\cos(18t))e^{({\rm i}t)}{\}},\\
\Gamma _{1}: {\{}z(t) &=& \xi+ (1.2+0.4\cos(4t))e^{(-{\rm i}t)}{\}},\qquad t: 0 \le t \le 2\pi,
\end{eqnarray*}
The test regions are  shown in Figure~\ref{Ex3Ahl} for two different values of $\xi$.  Given a first zero $a_0$ of the Ahlfors map, the exact second zero $a_1$ is also unknown for this region. We compute $a_{1n}$ from (\ref{kkk}) which is the approximate value of $a_1$. The acuracy is measured by computing $f(a_{1n})$ from (\ref{eq39}). The theoretical value of $f(a_1)$ is zero. See Table~\ref{tabEx3} for the numerical values of $a_{1n}$ and $f(a_{1n})$.

\end{example}

\begin{table}[h!]
\caption{Numerical values of $a_{1n}$ and $f(a_{1n})$ for Example~\ref{Ex3}.}\label{tabEx3}
\begin{center}
\begin{tabular}{ c c c c c}
\hline
\multicolumn{3}{c}{\rule{0ex}{3ex}$\xi = 0.6452-0.8655{\rm i}~\&~ a_0=-1.3088+1.8012{\rm i}$ } \\
\hline
$n$       &$a_{1n}$            &$f(a_{1n})$\\
\hline
$128$  &$2.0567 - 2.4889{\rm i}$&$3.0(-02)$\\
\hline
$256$   &$2.0625 - 2.5444{\rm i}$&$5.8(-05)$\\
\hline
$512$&$2.0624 - 2.5445{\rm i}$&$8.9(-08)$\\
\hline
$1024$ &$2.0624 - 2.5445{\rm i}$ &$9.8(-12)$\\
\hline
\multicolumn{3}{c}{\rule{0ex}{3ex} $\xi=-2.4516+2.3626{\rm i}~ \&~ a_0=2.2673-2.0351{\rm i}$} \\
\hline
\hline
$n$     &$a_{1n}$            &$f(a_{1n})$\\
\hline
$256$  &$-3.0647 + 3.2599{\rm i}$&$6.9(-04)$\\
\hline
$512$   &$-3.0648 + 3.2600{\rm i}$&$1.8(-07)$\\
\hline
$1024$  &$-3.0648 + 3.2600{\rm i}$  &$5.6(-12)$\\
\hline
\end{tabular}
\end{center}
\end{table}
\newpage
\section{Conclusion} \label{Concls}
In this paper, we have constructed a numerical method for finding second zero of the Ahlfors map of doubly connected regions. We derived two formulas for the derivative of the boundary correspondence function $\theta(t)$ of the Ahlfors map. These formulas were then used to find the second zero of the Ahlfors map for any smooth doubly connected regions.
Analytical method for computing the exact zeros of Ahlfors map for annulus region is presented in ~\cite{teg98} and ~\cite{teg} but the problem of finding  zeros for arbitrary doubly connected regions is the first time presented in this paper.  The numerical examples presented have illustrated that our method involving boundary integral equation has high accuracy.
\section*{Acknowledgment} \label{sc:int}
This work was supported in part by the Malaysian Ministry of Education (MOE) through the Research Management Centre (RMC), Universiti Teknologi Malaysia (FRGS Ref. No. PY/2014/04077). This support is gratefully acknowledged.

\end{document}